\documentclass{article}

\usepackage{proof}
\usepackage{latexsym}
\usepackage{amssymb}

\usepackage{color}

\newcommand{\bfthree}{\mbox{\boldmath$3$}}

\newcommand{\blem}{\begin{lemma}}
\newcommand{\elem}{\end{lemma}}
\newcommand{\bth}{\begin{theorem}}
\newcommand{\benu}{\begin{enumerate}}
\newcommand{\eenu}{\end{enumerate}}
\newcommand{\bdes}{\begin{description}}
\newcommand{\edes}{\end{description}}
\newcommand{\bdf}{\begin{definition}}
\newcommand{\edf}{\end{definition}}
\newcommand{\bcor}{\begin{cor}}
\newcommand{\ecor}{\end{cor}}
\newcommand{\bprp}{\begin{proposition}}
\newcommand{\eprp}{\end{proposition}}
\newcommand{\bmlem}{\begin{mlemma}}
\newcommand{\emlem}{\end{mlemma}}
\newcommand{\bclm}{\begin{claim}}
\newcommand{\eclm}{\end{claim}}
\newcommand{\bprf}{{\bf Proof}.\hspace{2mm}}

\newcommand{\eprf}{\hspace*{\fill} $\Box$}

\newcommand{\Integer}{\mathbb{Z}}
\newcommand{\bftwo}{\mbox{\boldmath$2$} }

\newcommand{\beqn}{\begin{equation}}
\newcommand{\eeqn}{\end{equation}}
\newcommand{\beqnarr}{\begin{eqnarray}}
\newcommand{\eeqnarr}{\end{eqnarray}}
\newcommand{\beqnarrs}{\begin{eqnarray*}}
\newcommand{\eeqnarrs}{\end{eqnarray*}}

\newcommand{\spand}{\,\&\,}

\newtheorem{theorem}{Theorem}[section]
\newtheorem{definition}[theorem]{Definition}
\newtheorem{proposition}[theorem]{Proposition}
\newtheorem{lemma}[theorem]{Lemma}
\newtheorem{cor}[theorem]{Corollary}
\newtheorem{mlemma}[theorem]{Main Lemma}
\newtheorem{claim}[theorem]{Claim}

\newcommand{\alp}{\alpha}

\newcommand{\del}{\delta}
\newcommand{\Del}{\Delta}
\newcommand{\ome}{\omega}

\newcommand{\bet}{\beta}
\newcommand{\gam}{\gamma}
\newcommand{\Gam}{\Gamma}

\newcommand{\Sig}{\Sigma}

\newcommand{\Tht}{\Theta}
\newcommand{\lam}{\lambda}
\newcommand{\Lam}{\Lambda}

\newcommand{\fal}{\forall}
\newcommand{\exi}{\exists}

\newcommand{\Rarw }{\Rightarrow}

\newcommand{\lrarw}{\leftrightarrow}
\newcommand{\Lrarw}{\Leftrightarrow}

\newcommand{\calE}{{\cal E}}

\newcommand{\calL}{{\cal L}}

\newcommand{\calP}{{\cal P}}

\newcommand{\calx}{{\cal X}}
\newcommand{\caly}{{\cal Y}}

\newcommand{\Natural}{\mathbb{N}}

\title{Cut-eliminability in second order logic calculi
}
\author{
Toshiyasu Arai
\\
Graduate School of Science,
Chiba University
\\
1-33, Yayoi-cho, Inage-ku,
Chiba, 263-8522, JAPAN
\\
tosarai@faculty.chiba-u.jp
}
\date{}
\begin{document}
\maketitle

\begin{abstract}
In this paper we propose a semantics in which the truth value of a formula is a pair of elements in a complete Boolean algebra.
Through the semantics we can unify largely two proofs of cut-eliminability (Hauptsatz) in classical second order logic calculus,
one is due to Takahashi-Prawitz and the other by Maehara.
\end{abstract}

\section{Takeuti's fundamental conjecture}
${\sf G}^{1}{\sf LC}$ defined in subsection \ref{subsec:GLC} below
is an impredicative sequent calculus with the $(cut)$ rule for the second order logic.
${\sf G}^{1}{\sf LC}^{cf}$ denotes the cut-free fragment of ${\sf G}^{1}{\sf LC}$, and
${\sf LK}={\sf G}^{0}{\sf LC}$ the first order fragment.
\\

\noindent
(Takeuti's fundamental conjecture for the second order calculus ${\sf G}^{1}{\sf LC}$\cite{Takeuti53})

$(cut)$ inferences are eliminable from proofs in ${\sf G}^{1}{\sf LC}$:
if ${\sf G}^{1}{\sf LC}$ proves a sequent, then it is provable without $(cut)$.
\\

It seems to me that G. Takeuti's intention in the conjecture is to reduce or paraphrase
the consistency problem of the second order arithmetic $\Integer_{2}=(\Pi^{1}_{\infty}\mbox{-CA})$
to a mathematical problem of cut-eliminability in the second order calculus ${\sf G}^{1}{\sf LC}$,
and the consistency of higher order arithmetic to the cut-eliminability in the higher order calculus {\sf GLC}.

Some partial results are obtained on the conjecture.
Takeuti \cite{FC5} shows a cut-elimination theorem for a fragment of ${\sf G}^{1}{\sf LC}$,
and one for a fragment of the higher order calculus {\sf GLC} in \cite{FC6},
both of which implies the 1-consistency of the subsystem $(\Pi^{1}_{1}\mbox{-CA})_{0}$ of the second order arithmetic,
the strongest one in the big five.
In \cite{SBL} a cut-elimination theorem for a fragment of ${\sf G}^{1}{\sf LC}$ is shown,
which implies the 1-consistency of the subsystem $(\Del^{1}_{2}\mbox{-CA+BI})$ of the second order arithmetic.
All of these proofs in \cite{FC5, FC6, SBL} are based on transfinite induction on computable notation systems of ordinals,
and hence are ordinal-theoretically informative ones.

Although no proof of the full conjecture has been obtained as Takeuti had expected,
the cut-eliminability holds for second order calculus.

\bth\label{th:Tait}{\rm \cite{Tait66}}\\
\[
{\sf G}^{1}{\sf LC}\vdash \Gam\Rarw\Del
\Rarw
{\sf G}^{1}{\sf LC}^{cf}\vdash \Gam\Rarw\Del
\]
\end{theorem}

Moreover the cut-eliminability holds for higher order calculus ${\sf GLC}$.

\bth\label{th:Takahashi}{\rm \cite{Takahashi, Prawitz}}\\
\[
{\sf G}{\sf LC}\vdash \Gam\Rarw\Del
\Rarw
{\sf G}{\sf LC}^{cf}\vdash \Gam\Rarw\Del
\]
\end{theorem}

In this paper let us focus on the second order calculus for simplicity,
and we propose a semantics in which the truth value of a formula is a pair of elements in a complete Boolean algebra.
Through the semantics we can unify largely two proofs of cut-eliminability (Hauptsatz) in classical second order logic calculus,
one is due to Takahashi-Prawitz and the other by Maehara.

In Section \ref{sect:val}
a soundness theorem \ref{cor:1} of ${\sf G}^{1}{\sf LC}$ is shown for semi valuations based on the semantics
with pairs of elements in a complete Boolean algebra.
Our proof of the theorem is essentially the same as in Takahashi\cite{Takahashi}, Prawitz\cite{Prawitz} and Maehara\cite{Maehara}.
In Section \ref{sect:Takahashi} Theorem \ref{th:Tait} is concluded.

In Section \ref{sect:Maehara} 
 a cBa $\mathbb{B}_{X}\subset\mathcal{P}(X)$ is introduced from a relation $M$ on an arbitrary set $X\neq\emptyset$.
 The construction of the cBa $\mathbb{B}_{X}$ is implicit in \cite{Maehara}.
Theorem \ref{th:Tait} is proved using a semi valuation defined from cut-free provability.
 
In Section \ref{sect:strength} the proof theoretic strength of cut-eliminability is calibrated.
It is well known that Theorem \ref{th:Tait} is equivalent to the 1-consistency of $\Integer_{2}$ over a weak arithmetic.
We sharpen it with respect to end sequents of proofs and fragments.
Finally some open problems are mentioned.

\subsection{Logic calculi}\label{subsec:GLC}

Let us recall second order sequent calculi briefly.
Details are found in \cite{Takeuti87}.

Logical connectives are $\lnot,\lor,\land,\exi,\fal$.
A second order language is obtained from a first order language by adding 
countably infinite $n$-ary variables $X^{n}_{i}\,(i\in\ome)$ for each $n=1,2,\ldots$.
For simplicity let us assume that our language contains no relation (predicate) symbol nor function symbol.
Formulas are quantified by second order quantifiers $\exi X^{n},\fal X^{n}$ as well as first order quantifiers $\exi x, \fal x$.
For a formula $G$ and a list $\vec{x}=(x_{1}\ldots,x_{n})$ of distinct variables, the expression $\lam \vec{x}.G$ 
is an $n$-ary \textit{abstract} or a \textit{term} of second order, and denoted by $T,\ldots$.
$Tm_{0}$ denotes the set of first order terms, and $Tm_{1}^{(n)}$ the set of $n$-ary abstracts $\lam \vec{x}.G(x_{1},\ldots,x_{n})$.

For formulas $F(X^{n})$, $F(\lam \vec{x}.G)$ denotes the formula up to renaming of bound variables,
obtained from $F$
by replacing each atomic formula $
X^{n}(t_{1},\ldots,t_{n})$ by $(\lam \vec{x}. G(x_{1},\ldots,x_{n}))(T)\equiv G(t_{1},\ldots,t_{n})$.

A finite set of formulas are said to be a \textit{cedent}, denoted $\Gam,\Del,\ldots$.
$\Gam,\Del:=\Gam\cup\Del$, $\Gam,A:=\Gam\cup\{A\}$.
A pair of cedents  $(\Gam,\Del)$ is denoted
$\Gam\Rarw\Del$, and called a \textit{sequent}.
$\Gam$ is said to be the \textit{antecedent}, $\Del$ \textit{succedent} of the sequent $\Gam\Rarw\Del$.

A sequent calculus ${\sf G}^{1}{\sf LC}$ is a logic calculus for the second order logic.
Its initial sequents are
\[
A,\Gam\Rarw\Del,A\, (A\mbox{:atomic})
\]
Inference rules are first order ones $(L\lnot), (R\lnot),(L\lor), (R\lor),(L\land),(R\land)$, 
$(L\exi^{0}),(R\exi^{0})$, $(L\fal^{0}),(R\fal^{0})$
\[
\begin{array}{cc}
\infer[(L\lnot)]{\lnot F,\Gam\Rarw \Del}
{
\lnot F,\Gam\Rarw \Del,F
}
&
\infer[(R\lnot)]{\Gam\Rarw\Del,\lnot F}
{
F,\Gam\Rarw \Del,\lnot F
}
\end{array}
\]
where $F$ is the \textit{minor formula},
and $\lnot F$ the \textit{major formula} of the inference rules $(L\lnot), (R\lnot)$.
\[
\begin{array}{cc}
\infer[(L\lor)]{F_{0}\lor F_{1},\Gam\Rarw \Del}
{
F_{0},F_{0}\lor F_{1},\Gam\Rarw \Del
&
F_{1},F_{0}\lor F_{1},\Gam\Rarw \Del
}
&
\infer[(R\lor)]{\Gam\Rarw\Del,F_{0}\lor F_{1}}
{
\Gam\Rarw \Del,F_{0}\lor F_{1}, F_{i}
}
\end{array}
\]
where $i\in\{0,1\}$, $F_{0},F_{1}$ are the \textit{minor formula},
and $F_{0}\lor F_{1}$ the \textit{major formula} of the inference rules $(L\lor), (R\lor)$.
\[
\begin{array}{cc}
\infer[(L\land)]{F_{0}\land F_{1},\Gam\Rarw \Del}
{
F_{i},F_{0}\land F_{1},\Gam\Rarw \Del
}
&
\infer[(R\land)]{\Gam\Rarw\Del,F_{0}\land F_{1}}
{
\Gam\Rarw \Del,F_{0}\land F_{1},F_{0}
&
\Gam\Rarw \Del,F_{0}\land F_{1},F_{1}
}
\end{array}
\]
where $i\in\{0,1\}$, $F_{0},F_{1}$ are the \textit{minor formula},
and $F_{0}\land F_{1}$ the \textit{major formula} of the inference rules $(L\land), (R\land)$.
\[
\begin{array}{cc}
\infer[(L\exi^{0})]{\exi x F(x),\Gam\Rarw \Del}
{
F(a),\exi x F(x),\Gam\Rarw \Del
}
&
\infer[(R\exi^{0})]{\Gam\Rarw\Del,\exi x F(x)}
{
\Gam\Rarw \Del,\exi x F(x), F(t)
}
\\
&
\\
\infer[(L\fal^{0})]{\fal x F(x),\Gam\Rarw \Del}
{
F(t),\fal x F(x),\Gam\Rarw \Del
}
&
\infer[(R\fal^{0})]{\Gam\Rarw\Del,\fal x F(x)}
{
\Gam\Rarw \Del,\fal x F(x),F(a)
}
\end{array}
\]
where in $(L\exi^{0}),(R\fal^{0})$, $a$ is an eigenvariable which does not occur in the lower sequent,
and $F(a)$ is the \textit{minor formula}.
In $(R\exi^{0}),(L\fal^{0})$, $t$ is a first order term, and $F(t)$ is the \textit{minor formula}.
$\exi x F(x)$ is the \textit{major formula} of the inference rules $(L\exi^{0}), (R\exi^{0})$,
and $\fal x F(x)$ the \textit{major formula} of the inference rules $(L\fal^{0}), (R\fal^{0})$.

The $(cut)$ inference
\[
\infer[(cut)]{\Gam,\Pi\Rarw\Del,\Tht}
{
\Gam\Rarw \Del,C
&
C,\Pi\Rarw\Tht
}
\]
There is no minor nor major formula of $(cut)$ inference.

Rules for second order quantifications $(L\exi^{1}),(R\exi^{1}),(L\fal^{1}),(R\fal^{1})$.
\[
\begin{array}{cc}
\infer[(L\exi^{1})]{\exi X F(X),\Gam\Rarw \Del}
{
F(Y),\exi X F(X),\Gam\Rarw \Del
}
&
\infer[(R\exi^{1})]{\Gam\Rarw\Del,\exi X^{n} F(X)}
{
\Gam\Rarw \Del,\exi X F(X), F(T)
}
\\
&
\\
\infer[(L\fal^{1})]{\fal X^{n} F(X),\Gam\Rarw \Del}
{
F(T),\fal X F(X),\Gam\Rarw \Del
}
&
\infer[(R\fal^{1})]{\Gam\Rarw\Del,\fal X F(X)}
{
\Gam\Rarw \Del,\fal X F(X), F(Y)
}
\end{array}
\]
where in $(L\exi^{1}),(R\fal^{1})$, $Y$ is an eigenvariable which does not occur in the lower sequent,
and $F(Y)$ is the \textit{minor formula}.
In $(R\exi^{1}),(L\fal^{1})$, $T$ is an $n$-ary second order term, and $F(T)$ is the \textit{minor formula}.
$\exi X F(X)$ is the \textit{major formula} of the inference rules $(L\exi^{1}), (R\exi^{1})$,
and $\fal X F(X)$ the \textit{major formula} of the inference rules $(L\fal^{1}), (R\fal^{1})$.

Since cedents here are finite sets of formulas, there are no explicit structural rules,
weakening (or thinning), contraction nor exchange in our sequent calculi.

\section{Valuations}\label{sect:val}
In this section let us propose a semantics in which the truth value of a formula is a pair of elements in a complete Boolean algebra,
and a soundness theorem \ref{cor:1} of ${\sf G}^{1}{\sf LC}$ is shown for semi valuations based on the semantics.

For a cBa (complete Boolean algebra) $\mathbb{B}$ let ${\tt D}\mathbb{B}$ denote the set of pairs
$(a,b)$ of elements $a,b\in\mathbb{B}$ such that $a\leq b$.
Here ${\tt D}$ stands for the axiom ${\tt D}\, : \, \Box A \to\Diamond A$ in the modal logic.
$-a$ denotes the complement of $a\in\mathbb{B}$.

\bdf
{\rm For a cBa $\mathbb{B}$ let
\[
{\tt D}\mathbb{B}:=\{(a,b)\in\mathbb{B}\times\mathbb{B}: a\leq b\}.
\]

Each ${\tt a}\in {\tt D}\mathbb{B}$ is written
${\tt a}=(\Box{\tt a},\Diamond{\tt a})$, where $\Box{\tt a}\leq\Diamond{\tt a}$.
For ${\tt a}, {\tt b}\in{\tt D}\mathbb{B}$ let
\[
\begin{array}{cc}
{\tt a}\leq{\tt b}  :\Lrarw  \Box{\tt a}\leq\Box{\tt b} \spand \Diamond{\tt a}\leq\Diamond{\tt b}
&
-{\tt a} := (-\Diamond{\tt a},-\Box{\tt a})
\\
{\tt a}\unlhd{\tt b}  :\Lrarw  \Box{\tt a}\leq\Box{\tt b} \spand \Diamond{\tt a}\geq\Diamond{\tt b}
&
\end{array}
\]
}
\edf
Then for $\{{\tt a}_{\lam}\}_{\lam}\subset{\tt D}\mathbb{B}$, the following hold.
\[
\begin{array}{cc}
\sup_{<}\{{\tt a}_{\lam}\}_{\lam} = (\sup_{\lam}\Box{\tt a}_{\lam}, \sup_{\lam}\Diamond{\tt a}_{\lam})
&
\inf_{<}\{{\tt a}_{\lam}\}_{\lam} = (\inf_{\lam}\Box{\tt a}_{\lam}, \inf_{\lam}\Diamond{\tt a}_{\lam})
\\
\sup_{\lhd}\{{\tt a}_{\lam}\}_{\lam} = (\sup_{\lam}\Box{\tt a}_{\lam}, \inf_{\lam}\Diamond{\tt a}_{\lam})
&
\inf_{\lhd}\{{\tt a}_{\lam}\}_{\lam} = (\inf_{\lam}\Box{\tt a}_{\lam}, \sup_{\lam}\Diamond{\tt a}_{\lam})
\end{array}
\]

Obviously ${\tt D}\mathbb{B}$ is a complete lattice under the order $\leq$ as well as under the order $\lhd$.
Note that
${\tt a}\unlhd{\tt b}\Lrarw (\Box{\tt a}\to \Box{\tt b})= (\Diamond{\tt b}\to \Diamond{\tt a})=1$,
where $(a\to b):=\sup\{-a,b\}$ for $a,b\in\mathbb{B}$.

For example for $\mathbb{B}=\bftwo=\{0,1\}$,
${\tt D}\bftwo$ is the set of three truth values $\bfthree:={\tt D}\bftwo=\{{\tt f},{\tt u},{\tt t}\}=\{(0,0),(0,1),(1,1)\}$,
where ${\tt f}<{\tt u}<{\tt t}$ and
${\tt u}\lhd{\tt f}, {\tt t}$.

\bprp\label{prp:1}{\rm (Monotonicity)}\\
For ${\tt a},{\tt b},{\tt a}_{\lam},{\tt b}_{\pi}\in {\tt D}\mathbb{B}$
\[
{\tt a}\unlhd{\tt b}\Rarw -{\tt a}\unlhd -{\tt b}
.\]
\[
\fal\pi\exi\lam({\tt a}_{\lam}\unlhd{\tt b}_{\pi}) \spand
\fal\lam\exi\pi({\tt a}_{\lam}\unlhd{\tt b}_{\pi})
\Rarw
\sup_{<}\{{\tt a}_{\lam}\}_{\lam}\unlhd\sup_{<}\{{\tt b}_{\pi}\}_{\pi}
\spand 
\inf_{<}\{{\tt a}_{\lam}\}_{\lam}\unlhd\inf_{<}\{{\tt b}_{\pi}\}_{\pi}
.\]
\eprp

\bdf
{\rm A ${\tt D}\mathbb{B}$\textit{-valued model} $\mathcal{M}$ is a pair $(D_{0},D_{1})$, where
$D_{0}\neq \emptyset$ is a non-empty set, $D_{1}=\bigcup_{n\geq 1}D_{1}^{(n)}$
and $D_{1}^{(n)}$ a non-empty set of functions $\alp:D_{0}^{n}\to{\tt D}\mathbb{B}$ for each $n=1,2,\ldots$

For each $\alp\in D_{1}^{(n)}$ introduce an $n$-ary relation constant $\bar{\alp}$, 
and each $t\in D_{0}$ is identified with the individual constant for $t$.
For formulas $A$ and $n$-ary abstracts $\lam \vec{x}.G(x_{1},\ldots,x_{n})$ with $t\in D_{0}$ and $\bar{\alp}$,
let us define recursively 
$\mathcal{M}(A)\in {\tt D}\mathbb{B}$, $\mathcal{M}(\lam \vec{x}.G(x_{1},\ldots,x_{n})): D_{0}^{n}\to{\tt D}\mathbb{B}$ as follows.
\benu
\item
$\mathcal{M}(\bar{\alp}(t_{1},\ldots,t_{n}))=\alp(t_{1},\ldots,t_{n})$ for $t_{1},\ldots,t_{n}\in D_{0}$.
$\mathcal{M}(\lnot F)=-\mathcal{M}(F)$.
\item
$\mathcal{M}(F_{0}\lor F_{1})=\sup_{<}\{\mathcal{M}(F_{0}), \mathcal{M}(F_{1})\}$. 
$\mathcal{M}(F_{0}\land F_{1})= \inf_{<}\{\mathcal{M}(F_{0}), \mathcal{M}(F_{1})\}$.
\item
$\mathcal{M}(\exi x F(x))= \sup_{<}\{\mathcal{M}(F(t)): t\in D_{0}\}$. 
\\
$\mathcal{M}(\fal x F(x))= \inf_{<}\{\mathcal{M}(F(t)): t\in D_{0}\}$.
\item
$\mathcal{M}(\exi X^{n} F(X))=\sup_{<}\{\mathcal{M}(F(\bar{\alp})): \alp\in D_{1}^{(n)}\}$.
\\
$\mathcal{M}(\fal X^{n} F(X))=\inf_{<}\{\mathcal{M}(F(\bar{\alp})): \alp\in D_{1}^{(n)}\}$.

\item
$\mathcal{M}(\lam \vec{x}.G(x_{1},\ldots,x_{n}))(t_{1},\ldots,t_{n})=\mathcal{M}(G(t_{1},\ldots,t_{n}))$.
\eenu
}
\edf

Intuitively $\mathcal{M}(A)=(a,b)$ means that
the degree of truth of $A$ is $a$, and one of non-falsity of $A$ is $b$.
When $\mathbb{B}=\bfthree$, $\Box\mathcal{M}(A)=1$ [$\Diamond\mathcal{M}(A)=1$] is related to 
the fact that $!A$ is valid [$?A$ is valid] in a three-valued structure for
Girard's three-valued logic with modal operators $!,?$ in \cite{GirardB}, resp.

For $\alp,\bet:D_{0}^{n}\to{\tt D}\mathbb{B}$ let
\[
\alp\unlhd\bet:\Lrarw
\fal \vec{t}\in D_{0}^{n}(\alp(\vec{t})\unlhd\bet(\vec{t}))
\]
and
\beqnarrs
&& \mathcal{M}\models 3CA  :\Lrarw 
\mbox{for each formula }G(x_{1},\ldots,x_{n},X^{k}) \mbox{ and each } \bet\in D_{1}^{(k)},
\\
&&
\mbox{ there exists an } \alp\in D_{1}^{(n)} 
\mbox{ such that }
\alp\unlhd\mathcal{M}(\lam \vec{x}.G(x_{1},\ldots,x_{n},\bar{\bet}))
\eeqnarrs

Recall that 
$Tm_{0}$ denotes the set of first order terms, and $Tm_{1}^{(n)}$ the set of $n$-ary abstracts $\lam \vec{x}.G(x_{1},\ldots,x_{n})$.

\bdf
{\rm Let $V$ be a map from the set of formulas $A$ to ${\tt D}\mathbb{B}$, $A\mapsto V(A)\in {\tt D}\mathbb{B}$.
$V$ is said to be
a} semi ${\tt D}\mathbb{B}$-valuation {\rm if it enjoys the following conditions:}
\benu
\item
$V(\lnot F)\unlhd -V(F)$.
\item
$V(F_{0}\lor F_{1})\unlhd\sup_{<}\{V(F_{0}), V(F_{1})\}$.
$V(F_{0}\land F_{1})\unlhd \inf_{<}\{V(F_{0}), V(F_{1})\}$.
\item
$V(\exi x F(x))\unlhd \sup_{<}\{V(F(t)): t\in Tm_{0}\}$.
\\
$V(\fal x F(x))\unlhd \inf_{<}\{V(F(t)): t\in Tm_{0}\}$.
\item
$V(\exi X^{n} F(X))\unlhd\sup_{<}\{V(F(T)): T\in Tm^{(n)}_{1}\}$.
\\
$V(\fal X^{n} F(X))\unlhd \inf_{<}\{V(F(T)): T\in Tm^{(n)}_{1}\}$.
\eenu
\edf

\bdf
{\rm
Let
\[
\mathbb{B}_{\Del}:=\{(a,a): a\in\mathbb{B}\}
.\]
A $\mathbb{B}$-valued model $\mathcal{N}$ is a pair $(D_{0},I)$ such that $D_{0}$ is a non-empty set,
$I=\bigcup I^{(n)}$ and $I^{(n)}$ is a non-empty set of functions $\calx :D_{0}^{n}\to\mathbb{B}_{\Del}$.

Let $\mathcal{N}(A):=\Box\mathcal{N}(A)=\Diamond\mathcal{N}(A)$ for any formula $A$, and
\beqnarrs
&& \mathcal{N}\models 2CA :\Lrarw
\\
&&
\mbox{for each } n\geq 1 \mbox{ and each formula } G(x_{1},\ldots,x_{n},X),
\\
&&
\mathcal{N}(\fal X\exi Y^{n}\fal x_{1},\ldots,x_{n}(Y(x_{1},\ldots,x_{n})\lrarw G(x_{1},\ldots,x_{n},X)))=1
\eeqnarrs
where $1$ denotes the largest element in $\mathbb{B}$.
}
\edf

\bprp\label{prp:semival}
Let $\mathbb{B}$ be a cBa.
\benu
\item\label{prp:semival1}
Suppose ${\sf G}^{1}{\sf LC}^{cf}\vdash\Gam\Rarw\Del$.
Then $\inf\{\Box V(A):A\in\Gam\}\leq\sup\{\Diamond V(B): B\in\Del\}$, i.e.,
$\Diamond V(\bigwedge\Gam\supset\bigvee\Del)=1$
for any semi ${\tt D}\mathbb{B}$-valuation $V$.

\item\label{prp:semival2}
Suppose ${\sf G}^{1}{\sf LC}\vdash\Gam\Rarw\Del$.
Then $\inf\{\mathcal{N}(A):A\in\Gam\}\leq\sup\{\mathcal{N}(B): B\in\Del\}$, i.e.,
$\mathcal{N}(\bigwedge\Gam\supset\bigvee\Del)=1$
for any $\mathbb{B}$-valued model $\mathcal{N}=(D_{0},I)$ with $\mathcal{N}\models 2CA$.
\eenu
\eprp

\blem\label{lem:2}{\rm (Cf. \cite{GirardB}.)}
\\
Let $V$ be a semi ${\tt D}\mathbb{B}$-valuation.
Define a ${\tt D}\mathbb{B}$-model $\mathcal{M}=(D_{0},D_{1})$ by
$D_{0}=Tm_{0}$ and $D_{1}^{(n)}=\{V(T)\in{}^{D_{0}^{n}}{\tt D}\mathbb{B}: T\in Tm^{(n)}_{1}\}$ with 
$V(\lam \vec{x}. G(x_{1},\ldots,x_{n}))(t_{1},\ldots,t_{n}):=V(G(t_{1},\ldots,t_{n}))$ for $t_{1},\ldots,t_{n}\in Tm_{0}$.
Then
for formula $F(X^{n})$, $T\in Tm^{(n)}_{1}$, and $\alp=V(T)$
\beqn\label{eq:1}
V(F(T))\unlhd\mathcal{M}(F(\bar{\alp}))
\eeqn
and
\beqn\label{eq:2}
\mathcal{M}\models 3CA
\eeqn
\elem
\bprf
(\ref{eq:1}):
This is seen by induction on formulas $F(X)$ using Proposition \ref{prp:1}.
For example consider the case $F(X)\equiv(\exi Y^{k}\, G(Y,X))$.
By the induction hypothesis we have 
$V(G(S,T))\unlhd\mathcal{M}(G(\bar{\bet},\bar{\alp}))$  for any $S\in Tm_{1}^{(k)}$ and $ \bet=V(S)$.
Then
$V(F(T))\unlhd\sup_{<}\{V(G(S,T)): S\in Tm_{1}^{(k)}\}\unlhd\sup_{<}\{\mathcal{M}(G(\bar{\bet},\bar{\alp})): \bet=V(S), S\in Tm_{1}^{(k)}\}=
\mathcal{M}(F(\bar{\alp}))$.
\\
(\ref{eq:2}):
Let $G(x_{1},\ldots,x_{n},Y^{k})$ be a formula.
For a $k$-ary abstract $T\in Tm_{1}^{(k)}$, let $\alp=V(T)$, $\bet=V(\lam \vec{x}. G(x_{1},\ldots,x_{n},T))\in D_{1}^{(n)}$.
From (\ref{eq:1}) we see for any $t_{1},\ldots,t_{n}\in Tm_{0}$ that
$\bet(t_{1},\ldots,t_{n})=V(\lam \vec{x}. G(x_{1},\ldots,x_{n},T))(t_{1},\ldots,t_{n})=V(G(t_{1},\ldots,t_{n},T))\unlhd\mathcal{M}(G(t_{1},\ldots,t_{n},\bar{\alp}))=\mathcal{M}(\lam \vec{x}. G(x_{1},\ldots,x_{n},\bar{\alp}))(t_{1},\ldots,t_{n})$.
\eprf

\blem\label{lem:3}{\rm (Cf. \cite{Takahashi, Prawitz, Maehara}.)}
\\
Let $\mathcal{M}=(D_{0},D_{1})$ be a ${\tt D}\mathbb{B}$-valued model such that $\mathcal{M}\models 3CA$.
Let for $\alp\in D_{1}^{(n)}$
\[
\begin{array}{cc}
I(\alp)  :=  \{\mathcal{X}\in{}^{D_{0}^{n}}\mathbb{B}_{\Del}: \alp\unlhd\mathcal{X}\}
&
I^{(n)}  :=  \bigcup\{I(\alp):\alp\in D_{1}^{(n)}\}
\end{array}
\]
Then for the
$\mathbb{B}$-valued model $\mathcal{N}=(D_{0},I)$, $\alp\in D_{1}^{(n)}$, $\calx\in I^{(n)}$ and 
formulas $F(X^{n})$, the following hold:
\beqn\label{eq:3}
\alp\unlhd\mathcal{X} \Rarw \mathcal{M}(F(\bar{\alp}))\unlhd\mathcal{N}(F(\overline{\mathcal{X}}))
\eeqn
and
\beqn\label{eq:4}
\mathcal{N}\models 2CA
\eeqn
\elem
\bprf
Note that $\alp\unlhd\mathcal{X} \Lrarw \fal t_{1},\ldots,t_{n}\in D_{0}[\Box\alp(t_{1},\ldots,t_{n})\leq\mathcal{X}(t_{1},\ldots,t_{n})\leq\Diamond\alp(t_{1},\ldots,t_{n})]$.
\\
(\ref{eq:3}):
This is seen by induction on formulas $F(X)$.
The case when $F(X)$ is an atomic formula $X(t_{1},\ldots,t_{n})$ is seen from the assumption $\alp\unlhd\mathcal{X}$.
Other cases follow from Proposition \ref{prp:1}.
For example consider the case $F(X)\equiv(\exi Y^{k}\, G(Y,X))$.
By the induction hypothesis we have 
$\mathcal{M}(G(\bar{\bet},\bar{\alp}))\unlhd\mathcal{N}(G(\overline{\caly},\overline{\calx}))$  
for any $\bet\in D_{1}^{(k)}$ and $\caly\in I^{(k)}$ with $\bet\unlhd\caly$.
On the other hand we have
$\fal \caly\in I^{(k)}\exi \bet\in D_{1}^{(k)}(\bet\unlhd\caly)$ and $\fal\bet\in D_{1}^{(k)}\exi \caly\in I^{(k)}(\bet\unlhd\caly)$
by the definition of $I^{(k)}$.
Hence Proposition \ref{prp:1} yields
$\mathcal{M}(F(\bar{\alp}))=\sup_{<}\{\mathcal{M}(G(\bar{\bet},\bar{\alp})): \bet\in D_{1}^{(k)}\}\unlhd
\sup_{<}\{\mathcal{N}(G(\overline{\caly},\overline{\calx})): \caly\in I^{(k)}\}=\mathcal{N}(F(\overline{\mathcal{X}}))$.
\\
(\ref{eq:4}):
For formulas $G(x_{1},\ldots,x_{n},X^{k})$ we need to show that
\\
$\mathcal{N}(\fal X\exi Y\fal \vec{x}(Y(x_{1},\ldots,x_{n})\lrarw G(x_{1},\ldots,x_{n},X)))=1$.
Let $\mathcal{X}\in I^{(k)}$, and
$D_{1}^{(k)}\ni\alp\unlhd\mathcal{X}$.
From $\mathcal{M}\models 3CA$ pick a
$\bet\in D_{1}^{(n)}$ such that $\bet\unlhd\mathcal{M}(\lam \vec{x}. G(x_{1},\ldots,x_{n},\bar{\alp}))$.
On the other hand we have
$\mathcal{M}(\lam \vec{x}. G(x_{1},\ldots,x_{n},\bar{\alp}))\unlhd\mathcal{N}(\lam \vec{x}. G(x_{1},\ldots,x_{n},\mathcal{X}))$
by (\ref{eq:3}).
Now let $\mathcal{Y}(t_{1},\ldots,t_{n})=\mathcal{N}(G(t_{1},\ldots,t_{n},\mathcal{X}))$.
Then $\bet\unlhd\mathcal{Y}$ and $\mathcal{Y}\in I^{(n)}$.
Therefore
$\mathcal{N}(\fal \vec{x}(\overline{\mathcal{Y}}(x_{1},\ldots,x_{n})\lrarw G(x_{1},\ldots,x_{n},\overline{\mathcal{X}})))=1$.
\eprf

\bth\label{cor:1}
Suppose ${\sf G}^{1}{\sf LC}\vdash\Gam\Rarw\Del$.
Then for any cBa $\mathbb{B}$,
$\inf\{\Box V(A):A\in\Gam\}\leq\sup\{\Diamond V(B): B\in\Del\}$, i.e.,
$\Diamond V(\bigwedge\Gam\supset\bigvee\Del)=1$
for any semi ${\tt D}\mathbb{B}$-valuation $V$.
\end{theorem}
\bprf
For a given semi ${\tt D}\mathbb{B}$-valuation $V$,
let $\mathcal{M}$ be the ${\tt D}\mathbb{B}$-model in Lemma \ref{lem:2}.
By (\ref{eq:1}) we see that
 $V(C)\unlhd\mathcal{M}(C)$ for formulas $C$.
Also $\mathcal{M}\models 3CA$ by (\ref{eq:2}).
Next let $\mathcal{N}$ be the $\mathbb{B}$-valued model in Lemma \ref{lem:3}.
(\ref{eq:3}) yields
 $\mathcal{M}(C)\unlhd\mathcal{N}(C)$.
Also $\mathcal{N}\models 2CA$ by  (\ref{eq:4}).
Now assume ${\sf G}^{1}{\sf LC}\vdash\Gam\Rarw\Del$.
We obtain
 $\inf\{\Box V(A):A\in\Gam\}\leq\inf\{\Box \mathcal{M}(A):A\in\Gam\}\leq\inf\{\mathcal{N}(A):A\in\Gam\}\leq\sup\{\mathcal{N}(B): B\in\Del\}\leq\sup\{\Diamond \mathcal{M}(B): B\in\Del\}\leq\sup\{\Diamond V(B): B\in\Del\}$
 by Proposition \ref{prp:semival}.\ref{prp:semival2}.
\eprf
\\

Although the intermediate step with ${\tt D}\mathbb{B}$-models in Lemma \ref{lem:2} due to J. Y. Girard is intuitively appealing,
it is dispensable.
The following Lemma \ref{lem:2+3} is seen as in Lemmas \ref{lem:2} and \ref{lem:3}.

\blem\label{lem:2+3}
Let $V$ be a semi ${\tt D}\mathbb{B}$-valuation.
Define a $\bftwo$-valued model $\mathcal{N}=(D_{0},I)$ with
$D_{0}=Tm_{0}$ as follows.

Let $T\equiv(\lam \vec{x}. G(x_{1},\ldots,x_{n}))\in Tm_{1}^{(n)}$.
Then $v(T)(t_{1},\ldots,t_{n}):=v(G(t_{1},\ldots,t_{n}))$
for $t_{1},\ldots,t_{n}\in Tm_{0}$, and
$I(T)=\{\mathcal{X}\in{}^{D_{0}^{n}}\bftwo: v(T)\unlhd\mathcal{X}\}$.
Let
$I^{(n)}=\bigcup\{I(T): T\in Tm_{1}^{(n)}\}$.

Then for $T\in Tm_{1}^{(n)}$, $\calx\in I^{(n)}$ and 
formulas $F(X^{n})$,
\beqn\label{eq:2+3l}
v(T)\unlhd\mathcal{X} \Rarw v(F(T))\unlhd\mathcal{N}(F(\overline{\mathcal{X}}))
\eeqn
and
\beqn\label{eq:2+3CA}
\mathcal{N}\models 2CA
\eeqn
\elem

\section{Semi valuation through proof search}\label{sect:Takahashi}
It is easy to conclude Theorem \ref{th:Tait} from Theorem \ref{cor:1} and the following Lemma \ref{lem:1}.
This is the proof by Takahashi\cite{Takahashi} and Prawitz\cite{Prawitz}.

\blem\label{lem:1}{\rm (Cf. \cite{Schuette60}.)}\\
Suppose ${\sf G}^{1}{\sf LC}^{cf}\not\vdash\Gam\Rarw\Del$.
Then there exists a semi $\bfthree$-valuation $V$ such that
$V(A)={\tt t}$ for $A\in\Gam$ and $V(B)={\tt f}$ for $B\in\Del$.
\elem
\bprf
By a canonical proof search,  we get an infinite binary tree of sequents supposing
 ${\sf G}^{1}{\sf LC}^{cf}\not\vdash\Gam\Rarw\Del$.
Pick an infinite path through the tree.
Let us define formulas occurring in antecedents of the path to be ${\tt t}$,
formulas occurring in succedents to be ${\tt f}$.
This results in a semi valuation $V(A)\in\bfthree={\tt D}\bftwo$ such that
$\fal A\in\Gam(V(A)={\tt t})$, $\fal B\in\Del(V(B)={\tt f})$.
\eprf
\\

\noindent
({\bf Proof} of Theorem \ref{th:Tait}, ver.1)
\\
Suppose ${\sf G}^{1}{\sf LC}^{cf}\not\vdash\Gam\Rarw\Del$.
By Lemma \ref{lem:1}, pick a semi $\bfthree$-valuation $V$ such that
$\fal A\in\Gam[V(A)={\tt t}]$, $\fal B\in\Del[V(B)={\tt f}]$.
Namely
$1=\inf\{\Box V(A):A\in\Gam\}\not\leq\sup\{\Diamond V(B): B\in\Del\}=0$.
Theorem \ref{cor:1} yieds
 ${\sf G}^{1}{\sf LC}\not\vdash\Gam\Rarw\Del$.
\eprf

\section{Semi valuation defined from cut-free provability}\label{sect:Maehara}
In this section following Maehara\cite{Maehara},
 a cBa $\mathbb{B}_{X}\subset\mathcal{P}(X)$ is first introduced from a relation $M$ on an arbitrary set $X\neq\emptyset$.
 $M$ is a symmetric relation such that if $(x,x)\in M$, then $(x,y)\in M$ for any $y\in X$.
 The construction of the cBa $\mathbb{B}_{X}$ is implicit in \cite{Maehara}.
 Second the Hauptsatz for ${\sf G}^{1}{\sf LC}$ is concluded using a semi valuation defined from cut-free
 provability as in \cite{Maehara}.
 
It seems to me that 
Maehara' s proof compares more straightforward with the proof in Section \ref{sect:Takahashi} due to Takahashi-Prawitz
in the sense that the latter proves the contraposition of the Hauptsatz.
The cost we have to pay is to elaborate a cBa from relations in Subsection \ref{subsec:cBa},
which gives an inspiration to researches in non-classical logics, e.g., cf.\,\cite{BJOno}.

\subsection{complete Boolean algebras induced from relations}\label{subsec:cBa}

Let $X\neq\emptyset$ be a non-empty set, and $M:X\ni x\mapsto M(x)\subset X$ a map.
Assume $M$ enjoys the following two conditions for any $x,y\in X$:
\beqnarr
x\in M(x) & \Lrarw & M(x)=X
\label{eq:cBa1}
\\
x\in M(y) & \Lrarw & y\in M(x)
\label{eq:cBa2}
\eeqnarr
Then let
\[
\mathbb{B}_{X}:=\{\alp\subset X: \alp=\bigcap\{M(x): \alp\subset M(x), x\in X\}\}
.\]
In the following we consider only subsets of the set $X$.
Let $\bigcap\emptyset:=X$.

\blem\label{eq:cBa7}
$\fal\alp\subset X[\bigcap\{\gam\in\mathbb{B}_{X}: \alp\subset\gam\}=\bigcap\{M(x): \alp\subset M(x)\}\in\mathbb{B}_{X}]$.
\elem
\bprf
Let $\bet=\bigcap\{\gam\in\mathbb{B}_{X}: \alp\subset\gam\}$, and $\del=\bigcap\{M(x): \alp\subset M(x)\}$.
First it is clear that $\alp\subset M(x) \Lrarw \del\subset M(x)$, and hence $\del\in\mathbb{B}_{X}$.

We show $\del\subset\bet$.
Assume $\alp\subset\gam\in\mathbb{B}_{X}$ and $\gam\subset M(x)$.
Then $\alp\subset M(x)$.
Hence $\del\subset M(x)$, and $\del\subset\bigcap\{M(x): \gam\subset M(x)\}=\gam$.
Thus $\del\subset\bet$.
\eprf

\bth\label{th:cBa}
$\mathbb{B}_{X}$ is a cBa with the following operations for $\alp,\bet\in\mathbb{B}_{X}$,
and $\{M(x):x\in X\}\subset\mathbb{B}_{X}$.
\benu
\item\label{th:cBa1}
$1=X$.
$0=\bigcap_{y\in X}M(y)=\{x\in X: x\in M(x)\}$.

\item\label{th:cBa2}
$\inf_{\lam}\alp_{\lam}=\bigcap_{\lam}\alp_{\lam}$ and $\alp\leq\bet\Lrarw \alp\subset\bet$.
\\
$\sup_{\lam}\alp_{\lam}=\bigcap\{\gam\in\mathbb{B}_{X}: \bigcup_{\lam}\alp_{\lam}\subset \gam\}$.

\item\label{th:cBa3}
complement 
$-\alp=\bigcap\{M(x): x\in\alp\}$.

\eenu
\end{theorem}
\bprf
It is clear that $\{M(x):x\in X\}\subset\mathbb{B}_{X}$.
\\
\ref{th:cBa}.\ref{th:cBa1}.
We show $\bigcap_{y\in X}M(y),X\in\mathbb{B}_{X}$.
$\bigcap_{y\in X}M(y)\in\mathbb{B}_{X}$ is obvious.
If $\exi x\in X(x\in M(x))$, then $X\in\mathbb{B}_{X}$ follows from (\ref{eq:cBa1}).
Otherwise $\bigcap\{M(x): X\subset M(x)\}=\bigcap\emptyset=X$.

Next we show
$\{x\in X: x\in M(x)\} \Rarw x\in\bigcap_{y\in X}M(y)$.
Assume $x\in M(x)$.
Then by (\ref{eq:cBa1}) $y\in X=M(x)$. (\ref{eq:cBa2}) yields $x\in M(y)$.
\\
\ref{th:cBa}.\ref{th:cBa2}.
Suppose $\{\alp_{\lam}\}_{\lam}\subset\mathbb{B}_{X}$.
Let $\bet=\bigcap\{M(y): \bigcap_{\lam}\alp_{\lam}\subset M(y)\}$.
We show $\bet\subset\alp_{\lam_{0}}=\bigcap\{M(x):\alp_{\lam_{0}}\subset M(x)\}$ for any $\lam_{0}$.
Let $\alp_{\lam_{0}}\subset M(x)$. Then $\bigcap_{\lam}\alp_{\lam}\subset M(x)$, and 
$\bet\subset M(x)$.
Hence $\bet\subset\alp_{\lam_{0}}$.
We obtain $\bet\subset\bigcap_{\lam}\alp_{\lam}$, and hence $\bet\in\mathbb{B}_{X}$.
Therefore $\inf_{\lam}\alp_{\lam}=\bigcap_{\lam}\alp_{\lam}$.
On the othe side we see
$\sup_{\lam}\alp_{\lam}=\bigcap\{\gam\in\mathbb{B}_{X}: \bigcup_{\lam}\alp_{\lam}\subset \gam\}$
from Lemma \ref{eq:cBa7}.
\\
\ref{th:cBa}.\ref{th:cBa3}.
Let $y\in\bigcap\{M(x):-\alp\subset M(x)\}$, and $x\in\alp$.
Then $-\alp\subset M(x)$, and $y\in M(x)$.
Hence $y\in\bigcap\{M(x):-\alp\subset M(x)\}\Rarw y\in -\alp$.
This means $-\alp\in\mathbb{B}_{X}$.

Next we show
\beqn\label{eq:cBa10}
\alp\subset M(x) \Rarw x\in -\alp
\eeqn
Assume $\alp\subset M(x)$ and $y\in\alp$. Then $y\in M(x)$, and $x\in M(y)$ by (\ref{eq:cBa2}).
Thus $x\in-\alp=\bigcap\{M(y): y\in\alp\}$.

Third we show $\alp\cap(-\alp)=0=\{x:x\in M(x)\}$.
Let $x\in \alp\cap(-\alp)$. Then $x\in M(x)$ by the definition of $-\alp$.
Conversely let $x\in M(x)$. Then $x\in \bigcap_{y}M(y)\subset\alp\subset X=M(x)$.
(\ref{eq:cBa10}) yields $x\in \alp\cap(-\alp)$.

Finally we show $\sup\{\alp,-\alp\}=X$.
Let $\alp,-\alp\subset\bet\in\mathbb{B}_{X}$.
If $\bet\subset M(x)$, then by (\ref{eq:cBa10}) we have $x\in -\alp\subset M(x)$.
(\ref{eq:cBa1}) yields $M(x)=X$.
Therefore $\bet=\bigcap\{M(x): \bet\subset M(x)\}=X$.
\eprf
\\

The complement $-M(y)$ of $M(y)$ is given in the following Proposition \ref{prp:cBa}.

\bprp\label{prp:cBa}
For $y\in X$,
let $m(y):=\bigcap\{M(x): y\in M(x)\}$.
Then

\beqn\label{eq:cBa8}
m(y)=-M(y)
\eeqn
\eprp
\bprf
(\ref{eq:cBa8}):
By Theorem \ref{th:cBa}.\ref{th:cBa3} and (\ref{eq:cBa2}) we have
$-M(y)=\bigcap\{M(x):x\in M(y)\}=\bigcap\{M(x):y\in M(x)\}=m(y)$.
\eprf

\subsection{semi valuation induced from relation}\label{subsec:Maehara}

In what follows let $X=S$ be the set of all sequents.

\bdf
{\rm For sequents $\Gam\Rarw\Del$
\[
M(\Gam\Rarw\Del):=\{(\Lam\Rarw\Tht)\in S : {\sf G}^{1}{\sf LC}^{cf}\vdash\Gam,\Lam\Rarw\Del,\Tht\}
.\]
}
\edf
It is clear that the map $S\ni x\mapsto M(x)\subset S$ enjoys (\ref{eq:cBa1}) and (\ref{eq:cBa2}).
(\ref{eq:cBa1}) follows from the contraction and weakening (thinning) rules,
while (\ref{eq:cBa2}) is seen from the exchange rule, all of these rules are implicit in our calculus ${\sf G}^{1}{\sf LC}^{cf}$.

Let $\mathbb{B}_{S}\subset \mathcal{P}(S)$ be the cBa induced by the map, cf.\,Theorem \ref{th:cBa}.
We have for sequents $x\in S$,
$x\in M(x) \Lrarw {\sf G}^{1}{\sf LC}^{cf}\vdash x$,
$0=\{x\in S: x\in M(x)\}=M(\Rarw)$ for the empty sequent $\Rarw$.

\bdf\label{df:VM}
{\rm For formulas $A$
\beqnarrs
\Diamond V(A) & := & M(\Rarw A)
\\
\Box V(A) & := & m(A\Rarw)=\bigcap\{M(\Gam\Rarw\Del): (A\Rarw)\in M(\Gam\Rarw\Del)\}
\eeqnarrs
}
\edf
By Theorem \ref{th:cBa} and (\ref{eq:cBa8}) in Proposition \ref{prp:cBa} we have $\Diamond V(A),\Box V(A)\in\mathbb{B}_{S}$.

\blem\label{lem:M}
$V$ is a semi ${\tt D}\mathbb{B}_{S}$-valuation.
\elem
\bprf
$\Box V(A)\subset \Diamond V(A)$ is seen from ${\sf G}^{1}{\sf LC}^{cf}\vdash A\Rarw A$.

The conditions of the $\Diamond$ are seen from the right rules.
\\
$\Diamond V(\exi X^{n} F(X))\supset\sup_{<}\{\Diamond V(F(T)): T\in Tm^{(n)}_{1}\}$:
From the rule $(R\exi^{2})$ we see that
$\Diamond V(\exi X F(X))=M(\Rarw \exi X^{n}F(X))\supset\bigcup_{T\in Tm^{(n)}_{1}}M(\Rarw F(T))$.
Hence $\mathbb{B}_{S}\ni \Diamond V(\exi X^{n} F(X))\supset
\bigcap\{\alp: \bigcup_{T\in Tm^{(n)}_{1}}M(\Rarw F(T))\subset \alp\}=\sup_{<}\{\Diamond V(F(T)): T\in Tm^{(n)}_{1}\}$.
\\
$\Diamond V(\fal X^{n} F(X))\supset\inf_{<}\{\Diamond V(F(T)): T\in Tm^{(n)}_{1}\}$:
$\Diamond V(\fal X F(X))=M(\Rarw \fal X F(X))\supset\bigcap_{T\in Tm^{(n)}_{1}}M(\Rarw F(T))=\inf_{<}\{\Diamond V(F(T)): T\in Tm^{(n)}_{1}\}$
is seen from the rule $(R\fal^{2})$.
\\
$\Diamond V(\lnot A)\supset \Diamond(-V(A))$:
By (\ref{eq:cBa8}) and the rule $(R\lnot)$, we obtain
$\Diamond(-V(A))=-\Box V(A)=-m(A\Rarw)=M(A\Rarw)\subset M(\Rarw \lnot A)=\Diamond V(\lnot A)$.

The conditions for $\Box$ are seen from the left rules using (\ref{eq:cBa8}) in Proposition \ref{prp:cBa}.
\\
$\Box V(\exi X^{n} F(X))\subset\sup_{<}\{\Box V(F(T)): T\in Tm^{(n)}_{1}\}$:
By (\ref{eq:cBa8}) it suffices to show that
$M(\exi X^{n} F(X)\Rarw)\supset \bigcap\{M(F(T)\Rarw): T\in Tm^{(n)}_{1}\}$,
which follows  from the rule $(L\exi^{2})$.
\\
$\Box V(\fal X^{n} F(X))\subset\inf_{<}\{\Box V(F(T)): T\in Tm^{(n)}_{1}\}$:
Again by (\ref{eq:cBa8}) it suffices to show that
$M(\fal X^{n} F(X)\Rarw)\supset\sup_{<}\{M(F(T)\Rarw): T\in Tm^{(n)}_{1}\}$, which follows from the rule $(L\fal^{2})$.
\\
$\Box V(\lnot A)\subset -\Diamond V(A)$: 
The rule $(L\lnot)$ yields $M(\lnot A\Rarw)\supset M(\Rarw A)$, from which and (\ref{eq:cBa8}) we obtain
$\Box V(\lnot A)=m(\lnot A\Rarw)=-M(\lnot A\Rarw)\subset -\Diamond V(A)$.
\eprf
\\

\noindent
({\bf Proof} of Theorem \ref{th:Tait}, ver.2)
\\
Suppose ${\sf G}^{1}{\sf LC}\vdash\Gam\Rarw\Del$.
From Theorem \ref{cor:1} and Lemma \ref{lem:M} we see that
$\bigcap\{\Box V(A):A\in\Gam\}\subset\sup\{\Diamond V(B): B\in\Del\}$
for the semi ${\tt D}\mathbb{B}_{S}$-valuation $V$ defined in Definition \ref{df:VM}.
Now we have $(\Gam\Rarw)\in \Box V(A)=m(A\Rarw)$ for any $A\in\Gam$ by weakening.
Hence $(\Gam\Rarw)\in \sup\{\Diamond V(B): B\in\Del\}=\bigcap\{M(x): \bigcup_{B\in\Del}M(\Rarw B)\subset M(x)\}$.
On the other hand we have $\bigcup_{B\in\Del}M(\Rarw B) \subset M(\Rarw\Del)$ by weakening.
Therefore $(\Gam\Rarw)\in M(\Rarw\Del)$, i.e., ${\sf G}^{1}{\sf LC}^{cf}\vdash\Gam\Rarw \Del$.
\eprf

\section{Proof-theoretic strengths}\label{sect:strength}

In the final section let us calibrate proof theoretic strengths of cut-eliminability.
For a class $\Phi$ of sequents $CE_{\Phi}({\sf G}^{1}{\sf LC})$ denotes the statement that any ${\sf G}^{1}{\sf LC}$-provable 
sequent in $\Phi$
is provable without the $(cut)$ rule.
When $\Phi$ is the set of all sequents, let $CE({\sf G}^{1}{\sf LC}):\Lrarw CE_{\Phi}({\sf G}^{1}{\sf LC})$.
$\mbox{{\rm I}}\Sig_{1}$ denotes the fragment of the first-order arithmetic 
in which the complete induction schema is restricted to $\Sig^{0}_{1}$-formulas in the language of first-order arithmetic.
Let $\Sig^{0}_{1}$ denote the set of $\Sig^{0}_{1}$-sequents in which no second-order quantifier occurs, and
first-order existential quantifier [first-order universal quantifier] occurs only positively [occurs only negatively], resp.
Then $1\mbox{{\rm -CON}}(\Integer_{2})$ denotes
the 1-consistency of the second order arithmetic $\Integer_{2}=(\Pi^{1}_{\infty}\mbox{-CA})$,
which says that every $\Integer_{2}$-provable $\Sig^{0}_{1}$-sequent is true.

\bth\label{th:strength}
\[
\mbox{{\rm I}}\Sig_{1}\vdash CE({\sf G}^{1}{\sf LC})\lrarw CE_{\Sig^{0}_{1}}({\sf G}^{1}{\sf LC})\lrarw 1\mbox{{\rm -CON}}(\Integer_{2})
.\]
\end{theorem}
\bprf
$(CE_{\Sig^{0}_{1}}({\sf G}^{1}{\sf LC})\to 1\mbox{{\rm -CON}}(\Integer_{2}))$.
This is shown in \cite{TakeutiRem} as follows.
Argue in $\mbox{{\rm I}}\Sig_{1}$.

Let $\calL^{2}$ denote the class of lower elementary recursive functions.
The class of functions contains the zero, successor, projection and modified subtraction functions 
and  is closed under composition and summation of functions. 
$\calL^2_*$ denotes the class of lower elementary recursive relations. 
Then it is easy, cf.\,\cite{Rose} to see that
the class $\calL^2_*$ is closed under boolean operations and bounded quantifications,
each function in $\calL^2$ is bounded by a polynomial, and
the truth definition of atomic formulas $R(x_{1},\ldots,x_n)$ for $R\in \calL^2_*$ is elementary recursive.

Suppose that $\Integer_{2}\vdash\exi x\, R$ for a $\Sig^{0}_{1}$-sentence $\exi x\, R$ with an $R\in\calL^{2}_{*}$.
In the $\Integer_{2}$-proof,
restrict each first-order quantifier $\fal x,\exi x$ to
$\fal x\in\Natural,\exi x\in\Natural$, where 
$\Natural(a):\equiv\fal X(X(0)\land\fal y(X(y)\supset X(Sy))\supset X(a))$ with the successor function $S$,
and $\exi x\in\Natural \,B:\lrarw(\exi x(\Natural(x)\land B))$, etc.
Let us denote the restriction of a formula $A$ by $A^{\Natural}$.
The comprehension axiom (CA) $\exi X\fal y(X(y)\lrarw G(y))$ follows from $(R\exi^{1})$.
Complete induction schema follows $\fal a\in \Natural\fal X(X(0)\land \fal y\in N(X(y)\supset X(Sy))\supset X(a))$.
We obtain a ${\sf G}^{1}{\sf LC}$-proof of a sequent $Eq^{\Natural},A_{0}^{\Natural}\Rarw \exi x\in\Natural\, R$
for an axiom $A_{0}$ of finitely many constants for functions in $\calL^{2}$
and the equality axiom $Eq:\Lrarw (\fal X\fal x,y(x=y\to (X(x)\lrarw X(y))))$.
$A_{0}$ is a universal formula $\fal x_{1},\ldots,x_{n} Q$ with a $Q\in\calL^{2}_{*}$.
Thus we obtain a ${\sf G}^{1}{\sf LC}$-proof of the sequent $Eq,A_{0}\Rarw \exi x\, R$.

Next let $\calE(a):\Lrarw(\fal X\fal y(a=y\to (X(a)\lrarw X(y))))$, and restrict
each first-order quantifier $\fal x,\exi x$ occurring in the  ${\sf G}^{1}{\sf LC}$-proof to
$\fal x\in \calE,\exi x\in \calE$.
Then we obtain a ${\sf G}^{1}{\sf LC}$-proof of the sequent $Eq^{\calE},A_{0}\Rarw \exi x\, R$,
where $Eq^{\calE}\Lrarw (\fal X\fal x,y\in\calE(x=y\to (X(x)\lrarw X(y))))$, which is provable.
Hence we obtain a ${\sf G}^{1}{\sf LC}$-proof of the sequent $A_{0}\Rarw \exi x\, R$.
Now by $CE_{\Sig^{0}_{1}}({\sf G}^{1}{\sf LC})$, i.e., the cut-eliminability from the proof with $\Sig^{0}_{1}$-end sequents, we get
${\sf G}^{1}{\sf LC}^{cf}\vdash A_{0}\Rarw \exi x\, R$, i.e., ${\sf LK}\vdash A_{0}\Rarw \exi x\, R$.
Then we see that $\exi x\, R$ is true.
\\

\noindent
$(1\mbox{{\rm -CON}}(\Integer_{2})\to CE({\sf G}^{1}{\sf LC}))$.
Although this is a folklore, cf.\,\cite{GirardB}, let us show it briefly.

It suffices to show in $\Integer_{2}$, the cut-eliminability from \textit{each} proof $P$ of a sequent $\Gam\Rarw\Del$ 
since the statement $CE({\sf G}^{1}{\sf LC})$ is a $\Pi^{0}_{2}$.
In what follows argue in $\Integer_{2}$, and consider the Takahashi-Prawitz' proof in Section \ref{sect:Takahashi} for simplicity.
First observe that Lemma \ref{lem:1} of the existence of a semi $\bfthree$-valuation $V$ is provable
\footnote{Maehara's proof in Section \ref{sect:Maehara} is formalizable in ${\sf ACA}_{0}$. $\alp\in{\tt D}\mathbb{B}_{S}$ is definable by an arithmetical formula.}
in ${\sf WKL}_{0}$, a fortiori in $\Integer_{2}$,
assuming that ${\sf G}^{1}{\sf LC}^{cf}\not\vdash\Gam\Rarw\Del$.

In Lemma \ref{lem:2+3} the satisfaction relation $\mathcal{N}\models F$ in the $\bftwo$-model $\mathcal{N}=(D_{0},I)$ is
second-order definable for each formula $F$.
Then for each formulas $F(X^{n})$ and $G(x_{1},\ldots,x_{n},X)$, we have
$v(T)\unlhd\mathcal{X} \Rarw v(F(T))\unlhd\mathcal{N}(F(\overline{\mathcal{X}}))$
and
\\
$\mathcal{N}(\fal X\exi Y^{n}\fal x_{1},\ldots,x_{n}(Y(x_{1},\ldots,x_{n})\lrarw G(x_{1},\ldots,x_{n},X)))=1$.
This suffices to evaluate the truth values of formulas occurring in the proof $P$,
and $\mathcal{N}(\Gam\Rarw\Del)=0$.
Hence $P$ is not a ${\sf G}^{1}{\sf LC}$-proof of the sequent $\Gam\Rarw\Del$.
A contradiction.
\eprf

\bprp\label{prp:CEsigma}
$\mbox{{\rm I}}\Sig_{1}\vdash CE_{\Sig^{0}_{1}}({\sf G}^{1}{\sf LC}) \to CE({\sf G}^{1}{\sf LC})$.
\eprp
\bprf
This follows from Theorem \ref{th:strength} indirectly.
Here is a direct proof.

P. P\"appinghaus\cite{Paeppinghaus} shows that
$\mbox{{\rm I}}\Sig_{1}\vdash CE_{\Pi^{1}}({\sf G}^{1}{\sf LC}) \to CE({\sf G}^{1}{\sf LC})$
by using cut-absorption and the joker translation,
where $\Pi^{1}$ denotes the set of sequents in which second-order universal quantifier [second-order existential quantifier]
occurs only positively [occurs only negatively], resp.
In what follows argue in $\mbox{{\rm I}}\Sig_{1}$.

Let $\Gam\Rarw\Del$ be a $\Pi^{1}$-sequent.
Erase each second-order quantifier $\fal X,\exi Y$ in the sequent to get a first-order sequent $\Gam_{0}\Rarw\Del_{0}$.
It is easy to see that
if ${\sf G}^{1}{\sf LC}\vdash\Gam\Rarw\Del$, then ${\sf G}^{1}{\sf LC}\vdash\Gam_{0}\Rarw\Del_{0}$, and
if ${\sf G}^{1}{\sf LC}^{cf}\vdash\Gam_{0}\Rarw\Del_{0}$, then ${\sf G}^{1}{\sf LC}^{cf}\vdash\Gam\Rarw\Del$.
Hence we obtain $CE_{\Pi^{0}}({\sf G}^{1}{\sf LC}) \to CE_{\Pi^{1}}({\sf G}^{1}{\sf LC})$
for the set $\Pi^{0}$ of first-order sequents.

Next let $H\in\Sig^{0}_{1}$ be an Herbrand normal form of the first-order formula $\bigwedge\Gam_{0}\supset\bigvee\Del_{0}$.
Then again it is easy to see that
if ${\sf G}^{1}{\sf LC}\vdash\Gam_{0}\Rarw\Del_{0}$, then ${\sf G}^{1}{\sf LC}\vdash\Rarw H$, and
if ${\sf LK}\vdash\Rarw H$, then ${\sf LK}^{cf}\vdash\Gam_{0}\Rarw\Del_{0}$.
Therefore $CE_{\Sig^{0}_{1}}({\sf G}^{1}{\sf LC}) \to CE_{\Pi^{0}}({\sf G}^{1}{\sf LC})$.
\eprf
\\

Let us mention a refinement for fragments.
$\Pi^{1}_{n}$ denotes the class of formulas $G\equiv(\fal X_{1}\exi X_{2}\cdots Q X_{n}\, A)$
with a first-order matrix $A$, and $Q=\fal$ when $n$ is odd, $Q=\exi$ else.
An abstract $T\equiv(\lam \vec{x}.G(x_{1},\ldots,x_{k}))$ is in $\Pi^{1}_{n}$ iff
$G\in\Pi^{1}_{n}$.
Then ${\sf G}^{1}{\sf LC}(\Pi^{1}_{n})$ denotes a fragment of the calculus ${\sf G}^{1}{\sf LC}$ in which
inference rules $(R\exi^{1}), (L\fal^{1})$ are restricted to $T\in\Pi^{1}_{n}$:
\[
\begin{array}{cc}
\infer[(R\exi^{1})]{\Gam\Rarw\Del,\exi X \, F(X)}
{
\Gam\Rarw\Del,\exi X\, F(X), F(T)
}
&
\infer[(L\fal^{1})]{\fal X\, F(X),\Gam\Rarw\Del}
{
F(T),\fal X\, F(X),\Gam\Rarw\Del
}
\end{array}
\]
An inspection to the proof of Theorem \ref{th:strength} shows the following.
Note that $\Natural(a)$ as well as $\calE(a)$ is a $\Pi^{1}_{1}$-formula without second-order free variable.

\bcor\label{cor:strength}
For each $n>0$
\[
\mbox{{\rm I}}\Sig_{1}\vdash CE({\sf G}^{1}{\sf LC}(\Pi^{1}_{n}))\lrarw 
CE_{\Sig^{0}_{1}}({\sf G}^{1}{\sf LC}(\Pi^{1}_{n}))\lrarw 1\mbox{{\rm -CON}}((\Pi^{1}_{n}\mbox{{\rm -CA}})_{0})
.\]
\ecor

Finally let us mention some open problems.
\\
{\bf Problem 1}.
What is the proof theoretic strength of the statement $CE_{\Pi^{0}_{1}}({\sf G}^{1}{\sf LC})$?
\\

$CE_{\Pi^{0}_{1}}({\sf G}^{1}{\sf LC})$ says that any  ${\sf G}^{1}{\sf LC}$-provable $\Pi^{0}_{1}$-sequent 
is provable without the $(cut)$ rule, where $\Pi^{0}_{1}$ denotes the dual class for $\Sig^{0}_{1}$.
Specifically does $\mbox{{\rm I}}\Sig_{1}$ prove $CE_{\Pi^{0}_{1}}({\sf G}^{1}{\sf LC})$?
\\

\noindent
To state the next problem we need first some definitions.

\bdf
{\rm
\benu

\item
An inference rule is said to be \textit{reducible}
if there is a minor formula $A$ of the inference rule such that
either the formula $A$ is in the antecedent and the sequent $\Rarw A$ is provable, or
$A$ is in the succedent and the sequent $A\Rarw$ is provable.

\item
A proof $P$ enjoys the \textit{pure variable condition} 
if in $P$, a free variable occurs in a sequent other than the end-sequent,
then it is an eigenvariable of an inference rule $J$
and the variable occurs only in the upper part of the inference rule $J$.

\item
A proof is said to be in \textit{irreducible} or in \textit{Mints' normal form}
 if it is cut-free, enjoys the pure variable condition
and contains no reducible inference rules.
\eenu
}
\edf

Mints' normal form theorem for a sequent calculus {\sf C} states that every {\sf C}-provable sequent has an irreducible
proof (with respect to {\sf C}).

\bth{\rm (\cite{Mints, AraiMints})}\\
Over $\mbox{{\rm I}}\Sig_{1}$,
Mints' normal form theorem for the first-order calculus {\sf LK} is equivalent to the 2-consistency of the first-order arithmetic {\sf PA}.
\end{theorem}
{\bf Problem 2}.
Does Mints' normal form theorem hold for ${\sf G}^{1}{\sf LC}$?
\\

It is easy to see that Mints' normal form theorem for ${\sf G}^{1}{\sf LC}$ implies the 2-consistency of 
the second-order arithmetic $\Integer_{2}$ as follows.
Assume that $\Integer_{2}\vdash\exi x\fal y\, R(x)$ for a false $\Sig^{0}_{2}$-sentence $\exi x\fal y \,R(x)$.
Let $Ind:\Lrarw(\fal a\,\Natural(a))$.
Then ${\sf G}^{1}{\sf LC}\vdash \Rarw\exi x(Eq\land Ind\land A_{0}\supset\fal y\, R(x))$ for a true $\Pi^{0}_{1}$-sentence $A_{0}$.
Pick an irreducible proof $P$ of the sequent $\Rarw\exi x(Eq\land Ind\land A_{0}\supset\fal y\, R(x))$ in ${\sf G}^{1}{\sf LC}$.
Then for a closed term $t$ the last inference must be a right rule $(R\exi^{0})$:
\[
\infer[(R\exi^{0})]{\Rarw\exi x(Eq\land Ind\land A_{0}\supset\fal y\, R(y))}
{
\Rarw\exi x(Eq\land Ind\land A_{0}\supset\fal y\, R(y)), Eq\land Ind\land A_{0}\supset\fal y\, R(t)
}
\]
From the $\Sig^{0}_{1}$-completeness, we see for the false $\Pi^{0}_{1}$-sentence $\fal y R(t)$,
that there exists a proof of the sequent
$Eq\land Ind\land A_{0}\supset\fal y\, R(t)\Rarw$ even in the weak fragment {\sf BC} of ${\sf G}^{1}{\sf LC}$ 
defined in p.166, \cite{Takeuti87}, in which
the abstracts $T$ in the inference rules $(R\exi^{1}),(L\fal^{1})$ are restricted to variables and predicate constants.
This means that $P$ is reducible. A contradiction.
\\

Cut-elimination by absorption in \cite{Paeppinghaus} is useless to prove the Mints' normal form theorem
since in
\[
\infer[(L\fal^{1})]{\fal X(X\supset X),\Gam\Rarw\Del}
{
 \infer[(L\supset)]{A\supset A,\Gam\Rarw\Del}
 {
\Gam\Rarw \Del,A
&
A,\Gam\Rarw\Del
}
}
\]
$\Rarw \fal X(X\supset X)$ as well as $\Rarw A\supset A$ is provable, and both inferences 
$(L\fal^{1})$ and $(L\supset)$ are reducible.

A proof of Mints' normal form theorem hold for {\sf LK} in \cite{AraiMints} runs as follows.
Assume that  a sequent $\Gam_{0}\Rarw\Del_{0}$ has no irreducible proof.
By a proof search,  we get an infinite binary tree of sequents, where we don't analyze, e.g., a succedent formula
$\exi x\, A(x)$ for a term $t$ when its instance $A(t)$ can be refuted, i.e., $A(t)\Rarw$ is provable.
\[
\infer{\Gam\Rarw\Del,\exi x\, A(x)}
{
\Gam\Rarw\Del,\exi x\, A(x),A(t)
&
\infer*{A(t)\Rarw}{}
}
\]
Pick an infinite path $\calP$ through the tree.
Let $\calP_{a}$ [$\calP_{s}$] denote the set of formulas occurring in an antecedent [occurring in a succedent]
of a sequent on the path $\calP$, resp.
Let atomic formulas in $\calP_{a}$ to be true, and atomic formulas in $\calP_{s}$ to be false.
From the truth values of atomic formulas define a first-order structure $\mathcal{M}$.
\footnote{Here we need $\ome$-times iterated jump operations. 
In a canonical proof search for cut-free provability in {\sf LK},
we obtain a valuation from an infinite path, which enjoys the Tarski's conditions without appealing iterated jump operations.
This is known as the Kreisel's trick.}
Then we see by induction on formulas $A$ that
if $A\in\calP_{a}$, then $\mathcal{M}\models A$, and
if $A\in\calP_{s}$, then $\mathcal{M}\not\models A$.
In the case of unanalyzed formula as above,
$\mathcal{M}\not\models A(t)$ follows from the soundness of the calculus {\sf LK} for any first-order structures $\mathcal{M}$.

An obstacle in extending this proof to ${\sf G}^{1}{\sf LC}$ lies in the fact that we need first prove
(\ref{eq:2+3l}), and then (\ref{eq:2+3CA}) follows from (\ref{eq:2+3l}) in the proof of Lemma \ref{lem:2+3}.
However in proving (\ref{eq:2+3l}) for an infinite path obtained from
a search tree with respect to the non-existence of irreducible proof,
we need the soundness of ${\sf G}^{1}{\sf LC}$ for $\bftwo$-models $\mathcal{N}$, but the soundness holds only if 
the model $\mathcal{N}$ enjoys the Comprehension axiom.
In other words we need  (\ref{eq:2+3CA}) before we prove (\ref{eq:2+3l}), and
we are in a circle.

\end{document}